\documentclass{article}
\usepackage{graphicx}
\usepackage{amsmath}
\usepackage{amsfonts}
\usepackage{amssymb}
\usepackage{ifthen}

\newtheorem{theorem}{Theorem}[section]
\newtheorem{lemma}[theorem]{Lemma}
\newtheorem{corollary}[theorem]{Corollary}
\newtheorem{remark}[theorem]{Remark}
\newtheorem{definition}[theorem]{Definition}
\newtheorem{example}[theorem]{Example}

\newcommand{\DEFINED}[1]{{\bf #1}}

\newcommand{\bis}{\prime\prime}

\newcommand{\Ii}{{\bf i}}

\newcommand{\IfAf}[1]{\Ii \frac{2\pi}{#1}}

\newcommand{\fphase}[2]{\exp{\left(\Ii \frac{2\pi}{#2} #1\right)}}

\newcommand{\FKP}{FKP}

\newcommand{\FFKP}{FFKP}

\newcommand{\FKPs}{\mathcal{F}}

\newcommand{\NofPART}[1]{p(#1)}

\newcommand{\NofEQcl}[1]{P_{#1}}

\newcommand{\EQclass}[2]{\left[ #2 \right]_{\RELofTYPE{#1}}}

\newcommand{\RELofTYPE}[1]{\stackrel{#1}{\simeq}}

\newcommand{\RELofEQUI}{\simeq}

\newcommand{\ROWtype}[2]{%
                          \ifthenelse{\equal{#2}{0}}
                             {1}
                             {\ifthenelse{\equal{#2}{1}}
                                 {1/(#1)}
                                 {1/({#1}^{#2})}%
                             }%
                        }

\newcommand{\NofROWS}[4]{{#1}^{#2}_{\ROWtype{#3}{#4}}}

\newcommand{\INTROix}[2]{{\tilde{n}}_{#2}}  

\newcommand{\Nt}[1]{\INTROix{a}{#1}}

\newcommand{\Si}[1]{\tilde{#1}}

\newcommand{\ITEMa}{a)}
\newcommand{\ITEMb}{b)}
\newcommand{\ITEMc}{c)}

\newcommand{\PROOFstart}{{\bf Proof}\\}
\newcommand{\PROOFend}{$\blacksquare$}

\begin{document}

\title{Permutation Equivalence Classes of Kronecker Products of Unitary Fourier Matrices}

\author{Wojciech Tadej\\
         {\small e-mail: wtadej@wp.pl}
       }
\date{\today}

\maketitle

\begin{abstract}
  Kronecker products of unitary Fourier matrices play important role
  in solving multilevel circulant systems by a multidimensional Fast
  Fourier Transform. They are also special cases of complex Hadamard
  (Zeilinger) matrices arising in many problems of mathematics and
  thoretical physics. The main result of the paper is splitting the
  set of all kronecker products of unitary Fourier matrices into
  permutation equivalence classes. The choice of the permutation
  equivalence to relate the products is motivated by the quantum
  information theory problem of constructing maximally entangled
  bases of finite dimensional quantum systems. Permutation
  inequivalent products can be used to construct inequivalent, in a
  certain sense, maximally entangled bases.   
\end{abstract}

%
%

\section{Introduction}

Kronecker products of unitary Fourier matrices, that is of matrices of
the form:
\begin{equation}
  [F_n]_{i,j} = \frac{1}{\sqrt{n}} e^{\IfAf{n} (i-1)(j-1)}\ \ \ \ \
  i,j \in \{1,2,\ldots,n\}
\end{equation}
play important role in solving multilevel circulant systems, that is
linear systems with matrices being linear combinations of kronecker products of
circulant matrices with the same structure. As Fourier matrices
diagonalize circulant matrices, their kronecker products diagonalize
such linear combinations, so repeated application of a multidimensional Fast Fourier
Transform enables one to solve effectively multilevel circulant
systems. Systems of this sort occure in problems of solving
partial differential equations with periodic (or partially periodic)
boundary conditions, image restoration problems (two-level
circulants), problems of approximating multilevel Toeplitz matrices
by multilevel circulants, which in turn arise in any context where one
encounters a 'shift invariant property'. Some applications of
multilevel circulant and Toeplitz matrices are considered in
\cite{Jin02}, while \cite{VanLoan92} exposes computational schemes for
a multidimensional FFT.

This work, however, originates from an other research area, namely
that of finding and classifying all real and complex Hadamard
(Zeilinger) matrices, that is unitary matrices with equal moduli of
their entries. Hadamard matrices find numerous applications in various
mathematical problems (see e.g. \cite{Haagerup96},
\cite{BjorkFroberg91}, \cite{BjorkSaffariS95}) and
theoretical  physics (see e.g. \cite{ReckZeilinger94}, \cite{JexStenholm95},
\cite{Werner01}, \cite{Ivanovic81}). The
question of finding all real and complex Hadamard matrices is still
open, see \cite{Dita04}, \cite{Dita04prepr}, \cite{Petrescu97}, \cite{Sloane}
for recent developments.
\bigskip

This article resolves the question when a kronecker product of Fourier
matrices can be transformed into another such product by left and
right multiplying it by unitary permutation (and diagonal) matrices.
The particular problem, which motivated the author to discriminate
between various kronecker products of Fourier
matrices using the above criterion, 
presented in \cite{Chhajlany03}, arises in the quantum information
theory and can be described, in terms of linear algebra, as follows.
\bigskip

Consider two $d$ dimensional quantum systems $\mathcal{A}$ and
$\mathcal{B}$, with ${\mathbf{C}}^d$ being the state space for each of
them, and ${\mathbf{C}}^d \otimes {\mathbf{C}}^d = {\mathbf{C}}^{d^2}$
the state space of the composite system $\mathcal{AB}$.  
Let $[a_j] = [a_0,a_1,\ldots,a_{d-1}]$ 
and $[b_k] = [b_0,b_1,\ldots,b_{d-1}]$ be orthonormal bases of ${\mathbf{C}}^d$,
chosen for $\mathcal{A}$ and $\mathcal{B}$ respectively. Then \mbox{$[a_j
\otimes b_k]$},\ \ \mbox{$j,k \in \{0,1,\ldots,(d-1)\}$} is an orthonormal basis
for the state space ${\mathbf{C}}^{d^2}$ of $\mathcal{AB}$, $\otimes$ denoting the
kronecker product here and onwards.

A new basis $[f^H_{j,k}]$ for $\mathcal{AB}$ is generated from $[a_j
\otimes b_k]$ by sending a state vector $a_j \otimes b_k$ through the
'local operation' unitary quantum gate $G_1$ defined on the
basis vectors by
\begin{equation}
  G_1 (a_j \otimes b_k) = 
       \left( \sum_{i=0}^{d-1} H_{i,j} a_j \right) \otimes b_k
\end{equation}
 where $H$ is a unitary
$d \times d$ matrix, and then sending the result through the 'nonlocal
operation' unitary gate $G_2$ defined by:
\begin{equation}
  G_2( a_j \otimes b_k) = a_j \otimes b_{(j-k) \mod d}
\end{equation}
That is, $[f^H_{j,k}]$ is defined by:
\begin{equation}
  f^H_{j,k} = G_2 G_1 (a_j \otimes b_k)
\end{equation}

The basis $[f^H_{j,k}]$ is said to be maximally entangled iff for any
\mbox{$j,k \in \{0,1,\ldots,(d-1)\}$} the sum of the $d$ disjoint $d \times
d$ diagonal blocks of the $d^2 \times d^2$ projection matrix
$(f^H_{j,k}) \cdot (f^H_{j,k})^{\star}$ (so called partial trace of this
matrix) is equal to $\frac{1}{d} I$. A unitary matrix $H$ must have
equal moduli of its entries for $[f^H_{j,k}]$ to be maximally
entangled, i.e. $H$ must be a complex Hadamard matrix.

Next a bases equivalence relation is introduced. Two bases
$[f^{H_1}_{j,k}]$ and $[f^{H_2}_{j,k}]$ are $\equiv$ equivalent iff
one can transform one of these bases into the other by a 'local' unitary operation,
that is there exist unitary $d \times d$ matrices $U_1$ and $U_2$,
permutations $\pi_1,\pi_2 : \{0,1,\ldots,(d-1)\} \rightarrow
\{0,1,\ldots,(d-1)\}$ and phases $\phi_{j,k}$,\ \ \mbox{$j,k \in
\{0,1,\ldots,(d-1)\}$} such that
\begin{equation}
  \label{eq_bases_equivalence}
  f^{H_2}_{j,k} = 
       e^{\Ii \phi_{j,k}} \cdot (U_1 \otimes U_2) \cdot f^{H_1}_{\pi_1(j),\pi_2(k)}
\end{equation}
The phase factor $e^{\Ii \phi_{j,k}}$ is allowed since multiplying a
state vector of a quantum system by such a factor (or any complex
number) produces a vector representing the same state.

It is shown in \cite{Chhajlany03} that if $H_1,\ H_2$ are kronecker products
of unitary Fourier matrices then $[f^{H_1}_{j,k}] \equiv
[f^{H_2}_{j,k}]$  is equivalent to
the existence of permutation matrices $P_r,\ P_c$  such that 
$H_2 = P_r \cdot H_1 \cdot P_c$. 
\bigskip

Section \ref{sec_fkp_P_equivalence} of this article provides
classification of kronecker 
products of unitary Fourier matrices with respect to
whether one can permute one such product into another. Section
\ref{sec_fkp_PD_equivalence} 
extends the validity of this classification to the case when
additionally left and right multiplying by unitary diagonal matrices
(called phasing onwards) is allowed when transforming the products.

%
%

\section{Kronecker products of Fourier matrices, row types and their basic properties}
\label{sec_FKPs_and_row_types}

By abuse of notation, we will write a Fourier matrix instead of a
unitary Fourier matrix. This is defined by:


\begin{definition}
\label{def_Fourier_matrix}

  A \DEFINED{Fourier matrix} of size $n$ is an $n \times n$ matrix defined by:
  \begin{equation}
  \label{form_Fourier_matrix}
  [F_n]_{i,j} = \frac{1}{\sqrt{n}} e^{\IfAf{n} (i-1)(j-1)}\ \ \ \ \
  i,j \in \{1,2,\ldots,n\}
  \end{equation}

\end{definition}
\medskip

\begin{definition}
\label{def_FKP}

  A \DEFINED{Fourier Kronecker Product} (\FKP, \FKP\ product, \FKP\ subproduct) is a
  matrix $F$ such that
  \begin{equation}
    F = F_{N_1} \otimes F_{N_2} \otimes  \ldots  \otimes F_{N_r}
  \end{equation}
  where $r, N_1, N_2, \ldots ,N_r$ are some natural numbers and $F_k$
  is the Fourier matrix of size $k$.

\end{definition}
\medskip


\begin{definition}
\label{def_FFKP}

  A \DEFINED{Factored Fourier Kronecker Product} (\FFKP, \FFKP\ product, \FFKP\
  subproduct) is a matrix $F$ such that:
  \begin{equation}
    F = F_{P_1} \otimes F_{P_2} \otimes  \ldots  \otimes F_{P_r}
  \end{equation}
  where $r$ is a natural number, $P_1,P_2, \ldots ,P_r$ are natural
  powers of prime numbers, not necessarily distinct ones.

\end{definition}
\medskip


\begin{definition} 
\label{def_pFFKP}

  A \DEFINED{pure \FFKP} (p. \FFKP ) is an \FFKP\ (see Definition
  \ref{def_FFKP})  of the form, where $a$ is a prime number:
  \begin{equation}
    F = F_{a^{k_m}} \otimes F_{a^{k_{m-1}}} \otimes  \ldots  \otimes  F_{a^{k_1}}
  \end{equation}
  It is called \DEFINED{ordered} (p.o. \FFKP ) if $k_m \geq k_{m-1} \geq  \ldots  \geq k_1$.

\end{definition}
\medskip


\begin{definition}
\label{def_row_types}

  A \DEFINED{row of type} $\mathbf{\ROWtype{n}{1}}$, $n$ a natural number, of an \FKP\
  product of size $N$ (see Definition \ref{def_FKP}), is a row containing the elements:
  \begin{equation}
    \frac{1}{\sqrt{N}}\fphase{0}{n}, 
    \frac{1}{\sqrt{N}}\fphase{1}{n}, 
    \ldots, 
    \frac{1}{\sqrt{N}}\fphase{(n-1)}{n}
  \end{equation}
  each of which occures the same number of times in this row.

\end{definition}
\medskip


\begin{lemma}
\label{lem_four_rows_count}

  A Fourier matrix $F_{a^m}$, where $a$ is a prime number and $m$ is
  natural, contains only rows of types $\ROWtype{a}{0}, \ROWtype{a}{1}, \ROWtype{a}{2}, \ldots,
  \ROWtype{a}{m}$. The number of rows of each of those types is, respectively,
  $1,\ (a-1),\ (a-1)a,\ \ldots,\ (a-1)a^{m-1}$. That is, there are
  $(a-1)a^{r-1}$ rows of type $\ROWtype{a}{r}$.

\end{lemma}


\PROOFstart
Recall $F_n$ is given by $[F_n]_{i,j} = \frac{1}{\sqrt{n}} e^{\IfAf{n}
  (i-1)(j-1)}$. We will further concentrate on phases
$\frac{(i-1)(j-1)}{n} 2\pi$, more precisely on corresponding phase
fractions $\frac{(i-1)(j-1)}{n}$.


As the phase fraction of $[F_n]_{i,2}$ is $(i-1)/n$, for $1 < i \leq
n$ its proper ($<1$) fraction part is nonzero, so there is only one, the
first, row of $F_{a^m}$ containing full angle phases, that is a row of
type $\ROWtype{a}{0}$.


Let $1 \leq k \leq m$ and let $x \in \{0,1,\ldots,(a^k-1)\}$ be
relatively prime to $a^k$. We will show that the $(xa^{m-k} + 1)$-th
row of $F_{a^m}$ is a row of type $\ROWtype{a}{k}$.

The phase fraction sequence in the $(xa^{m-k} + 1)$-th row is the
sequence:
\begin{eqnarray*}
 \lefteqn{\left( 
           \frac{xa^{m-k}}{a^m} \cdot (j-1)\ :\ \ j = 1,2,\ldots,a^m
          \right)=}  \\
 & & \left( 
      \frac{x}{a^k} 0,\ \frac{x}{a^k} 1,\ \ldots,\ \frac{x}{a^k} (a^k-1),\ 
      \frac{x}{a^k} (a^k+0),\ \frac{x}{a^k} (a^k+1),\ \ldots,\ \frac{x}{a^k} (2a^k-1),\ 
      \ldots \right)
\end{eqnarray*}
with $a^{m-k}$ subsequences:
\begin{equation}
\label{form_rth_subseq_phasefrac}
  \frac{x}{a^k} (ra^k+0),\ \frac{x}{a^k} (ra^k+1),\ \ldots,\ \frac{x}{a^k} (ra^k+(a^k-1)) 
  \ \ \ \ \ \ \  r = 0,1,\ldots,(a^{m-k}-1)
\end{equation}
The proper fraction parts in the $r$-th subsequence are all different
and are equal to $\frac{0}{a^k},\ \frac{1}{a^k},\ \ldots,\
\frac{a^k-1}{a^k}$. If it were not so, we would have the equality:
\begin{equation}
  x(ra^k + p) \stackrel{\mod a^k}{=} x(ra^k + q)
\end{equation}
for some $p,q \in \{0,1,\ldots,(a^k-1)\},\ \ q > p$. But the above
equality implies that $a^k$ divides $x(q-p)$, that is $a^k$ divides
$(q-p)$ which is impossible. (Obviously, the proper fraction parts of
$p$-th's elements of (\ref{form_rth_subseq_phasefrac}) for a given $p$
and $r = 0,1,\ldots,(a^{m-k}-1)$ are all equal.) Thus the $(xa^{m-k} +
1)$-th row of $F_{a^m}$ is a row of type $\ROWtype{a}{k}$.


The number of the above considered rows is the number of $x$'s
belonging to $\{0,1,\ldots,(a^k-1)\}$ and relatively prime to
$a^k$. These are all the elements of $\{0,1,\ldots,(a^k-1)\}$ except
for multiples of $a$: $0a,\ 1a,\ 2a,\ \ldots,\ (a^{k-1} - 1)a$, that
is $a^k - a^{k-1} = a^{k-1}(a-1)$ numbers in total.

Since for a given $k$ the appropriate $x$'s index different rows, and for a
different $k$ we have other rows (of other type) indexed by
corresponding $x$'s, all $a^m$ rows of $F_{a^m}$ are indexed in this
way:
\begin{description}
  \item[] $1$ row of type $\ROWtype{a}{0}$
  \item[] $(a-1)$ rows of type $\ROWtype{a}{1}$\ \ \ \ 
          ($\frac{x}{a} \in
          \{\frac{1}{a},\frac{2}{a},\ldots,\frac{a-1}{a}\}$)
  \item[] $(a-1)a$ rows of type $\ROWtype{a}{2}$
  \item[] $\ldots$
  \item[] $(a-1)a^{m-1}$ rows of type $\ROWtype{a}{m}$
\end{description}
\PROOFend
\bigskip


\begin{lemma}
\label{lem_row_kronp_type}

  The kronecker product of two rows of types $\ROWtype{a}{k},\
  \ROWtype{a}{l}$ is a row of type $\ROWtype{a}{\max(k,l)}$.

\end{lemma}


\PROOFstart
Let $k \geq l$. A row of type $\ROWtype{a}{k}$ contains elements with
phase fractions $0/a^k,\ 1/a^k,\ \ldots,\ (a^k-1)/a^k$ multiplied by
$2\pi$, each phase occuring $p$ times. Similarly, a row of type
$\ROWtype{a}{l}$ has phase fractions $0/a^l,\ 1/a^l,\ \ldots,\
(a^l-1)/a^l$, each occuring $q$ times.

The kronecker product of such rows contains phase fractions $x/a^k +
y/a^l$, where $x \in \{0,1,\ldots,(a^k-1)\}$ and $y \in
\{0,1,\ldots,(a^l-1)\}$, each occuring $p \cdot q$ times.

Note that for a given $y$ the phase fractions $x/a^k + y/a^l$ for \mbox{$x
\in \{0,1,\ldots,(a^k-1)\}$} can be simplified (as multiplied by $2\pi$
they later make phases of the row elements) to the proper ($<1$) phase
fractions $0/a^k,\ 1/a^k,\ \ldots,\ (a^k-1)/a^k$, each occuring $p
\cdot q$ times, as the proper fraction parts of $x/a^k + y/a^l,\ \ \ x
\in \{0,1,\ldots,(a^k-1)\}$ are all different. Were it not so, we
would have for some $x_1 < x_2$:
\begin{equation}
  x_1  +  y a^{k-l}   \stackrel{\mod a^k}{=}   x_2  +  y a^{k-l}
\end{equation}
which implies that $x_2 - x_1$ is divided by $a^k$, which is
impossible as $x_1,x_2 \in \{0,1,\ldots,(a^k-1)\}$.

Thus, for each $y \in \{0,1,\ldots,(a^l-1)\}$ we have indicated $pq$
occurences of each of the phases 
$\frac{0}{a^k} 2\pi,\ \frac{1}{a^k} 2\pi,\ \ldots,\ \frac{a^k-1}{a^k}
2\pi$. Then the total number of occurences of each such phase in the
kronecker product of the considered rows  is $a^l \cdot pq$, so the
product row is a row of type $\ROWtype{a}{k}$, where $k \geq l$ is
assumed.\PROOFend
\bigskip


\begin{lemma}
\label{lem_rel_prime_type_row_kronp_type}

  The kronecker product of two rows of types $\ROWtype{p}{1},\
  \ROWtype{q}{1}$, where $p$ and $q$ are relatively prime natural numbers
  is a row of type $\ROWtype{pq}{1}$.

\end{lemma}


\PROOFstart
A row of type $\ROWtype{p}{1}$ contains the phase fractions $0/p,\
1/p,\ \ldots,\ (p-1)/p$ multiplied by $2\pi$, each phase occuring $m$
times. Similarly, a row of type $\ROWtype{q}{1}$ contains the phase
fractions $0/q,\ 1/q,\ \ldots,\ (q-1)/q$, each fraction occuring $n$
times. 

The kronecker product of such rows contains phase fractions $x/p +
y/q$, where $x \in \{0,1,\ldots,(p-1)\}$ and $y \in
\{0,1,\ldots,(q-1)\}$, each fraction occuring $m \cdot n$ times. We
will show that all those $p \cdot q$ fractions can be reduced (as
multiplied by $2\pi$ they later make phases of the row elements) to
their proper ($<1$) fraction parts, which are all distinct and equal to
$0/(pq),\ 1/(pq),\ \ldots,\ (pq-1)/(pq)$.

Assume that for some pairs $x_1,\ y_1$ and $x_2,\ y_2$ the
corresponding phase fractions $x_1/p + y_1/q$ and  $x_2/p + y_2/q$
have their proper fraction parts equal. That is
\begin{equation}
  x_1 q  +  y_1 p   \stackrel{\mod pq}{=}    x_2 q  +  y_2 p
\end{equation}
Then $pq$ divides 
$(x_2 q  +  y_2 p)  -  (x_1 q  +  y_1 p)  =  q(x_2 - x_1)  +  p(y_2 - y_1)$,
and as $p$ and $q$ are relatively prime, $p$ divides $(x_2 - x_1)$ and
$q$ divides $(y_2 - y_1)$. Since $x_1,x_2 \in
\{0,1,\ldots,(p-1)\}$ and $y_1,y_2 \in \{0,1,\ldots,(q-1)\}$, we must
have $x_1 = x_2$ and $y_1 = y_2$.

We have thus shown that the kronecker product of rows of types
$\ROWtype{p}{1}$ and $\ROWtype{q}{1}$, where $p,\ q$ relatively prime,
contains phase fractions\\ $0/(pq),\ 1/(pq),\ \ldots,\ \mbox{$(pq-1)/(pq)$}$,
each fraction occuring $m \cdot n$ times, so the product row is a row
of type $\ROWtype{pq}{1}$.\PROOFend
\bigskip

Using the above lemmas and the fact that rows of a kronecker product are kronecker products of rows of
kronecker product factors, one can easily prove the two lemmas
below. See also the definition and remarks that follow.
\medskip


\begin{lemma}
\label{lem_poFFKP_row_types}

  A pure ordered \FFKP\ (see Definition \ref{def_pFFKP}) product of the form:
  \begin{equation}
    F = F_{a^{k_m}} \otimes F_{a^{k_{m-1}}} \otimes  \ldots  \otimes  F_{a^{k_1}}
  \end{equation}
  contains only rows of types $\ROWtype{a}{0},\ROWtype{a}{1}, \ROWtype{a}{2}, \ldots,
  \ROWtype{a}{k_m}$

\end{lemma}
\medskip


\begin{lemma}
\label{lem_types_in_two_pFFKPs_kprod}

  An \FFKP\ product $F$ being a kronecker product of two pure ordered
  \FFKP's (see Definitions \ref{def_FFKP},\ref{def_pFFKP}), $a \neq b$ prime numbers:
  \begin{equation}
    F   = F' \otimes F'',\ \  
    F'  = F_{a^{k_m}} \otimes  \ldots  \otimes  F_{a^{k_1}},\ \ 
    F'' = F_{b^{l_n}} \otimes  \ldots  \otimes  F_{b^{l_1}}
  \end{equation}
  contains only rows of types $\ROWtype{a^k b^l}{1}$ where $k
  \in \{0,1,\ldots,k_m\}$,\ \  $l \in \{0,1,\ldots,l_n\}$.

\end{lemma}
\medskip

More generally, we have:
\medskip


\begin{lemma}

  An \FFKP\ product $F$ being a kronecker product of $r$ pure ordered
  \FFKP's for $r$ different primes $a_1,a_2,\ldots,a_r$:
  \begin{equation}
    F = F^{(1)} \otimes F^{(2)} \otimes \ldots \otimes F^{(r)}
  \end{equation}
  where
  \begin{equation}
    F^{(k)} = F_{a_{k}^{p_{m^{(k)}}^{(k)}}} \otimes \ldots \otimes
              F_{a_{k}^{p_{1}^{(k)}}}
  \end{equation}
  contains only rows of types $\ROWtype{a_1^{l_1} \ldots
  a_r^{l_r}}{1}$, where $l_k \in \{0,1,\ldots,
  p_{m^{(k)}}^{(k)} \}$.

\end{lemma}
\bigskip


\begin{definition}
\label{def_kron_multiindex}
  A \DEFINED{kronecker multiindex} $i_1,\ldots,i_r;j_1,\ldots,j_r$, where
  $i_k,j_k \in \{1,\ldots,n_k\},\ \ n_k \in \mathbf{N}$ into a square
  matrix $A$ of size $n = n_1 \cdot n_2 \cdot \ldots \cdot n_r$
  (corresponding to the potential kronecker product structure of $A =
  A_1 \otimes A_2 \otimes \ldots \otimes A_r$, where $A_k$ of size
  $n_k$) indicates the $i,j$-th element of $A$, where:
  \begin{eqnarray}
    i & = & (i_1     - 1)\prod_{k=2}^{r} n_k  +
            (i_2     - 1)\prod_{k=3}^{r} n_k  + \ldots +
            (i_{r-1} - 1) n_r                 + i_r      \nonumber \\
    j & = & (j_1     - 1)\prod_{k=2}^{r} n_k  +
            (j_2     - 1)\prod_{k=3}^{r} n_k  + \ldots +
            (j_{r-1} - 1) n_r                 + j_r      \nonumber
  \end{eqnarray}
  We write $i_1,\ldots,i_r$ or $j_1,\ldots,j_r$  to kronecker
  mutliindex a row or column, respectively.

\end{definition}


\begin{remark}
  Let $H = H_1 \otimes H_2 \otimes \ldots \otimes H_r$ be a kronecker
  product of matrices $H_k$ of size $n_k$. Then the
  $i_1,\ldots,i_r;j_1,\ldots,j_r$-th element of $H$ is given by:
  \begin{equation}
    H_{i_1,\ldots,i_r;j_1,\ldots,j_r} = 
      [H_1]_{i_1,j_1} \cdot 
      [H_2]_{i_2,j_2} \cdot \ldots \cdot
      [H_r]_{i_r,j_r}  
  \end{equation}
\end{remark}


\begin{remark}
\label{rem_reorder_kronp}

   Reordering factors of a kronecker product is equivalent
   to left and right multiplying this product by appropriate
   permutation matrices (thus reordering doesn't change the row types
   present in an \FKP).

   These permutations are defined in the following way. 
   Suppose we reorder 
   $H = H_1 \otimes \ldots \otimes H_r$, $H_k$ of size $n_k$, 
   into 
   $H' = H_{\sigma(1)} \otimes \ldots \otimes H_{\sigma(r)}$, 
   where $\sigma$ is a permutation. 
   Then the $i,j$-th element of $H$, kronecker multiindexed also by
   $i_1,\ldots,i_r;j_1,\ldots,j_r$ (see Definition \ref{def_kron_multiindex}),
   where $i_k,j_k \in \{1,\ldots,n_k\}$, is moved by the appropriate
   permutation matrices into the
   $i_{\sigma(1)},\ldots,i_{\sigma(r)};j_{\sigma(1)},\ldots,j_{\sigma(r)}$-th 
   position within $H'$.

\end{remark}
\bigskip

%
%

\section{Permutation equivalence of kronecker \mbox{products} of Fourier matrices}
\label{sec_fkp_P_equivalence}

In this section we will basically split the set of all possible
kronecker products of Fourier matrices into equivalence classes with
respect to the permutation equivalence relation, denoted by
$\RELofTYPE{P}$ and defined below. This is an equivalence relation.


\begin{definition}
\label{def_permutation_equivalence}

  Two square matrices $A$ and $B$ of the same size are \DEFINED{permutation
    equivalent}, i.e. $A \RELofTYPE{P} B$,
  if and only if there exist permutation matrices $P_r,\ P_c$ of the
  proper size such that
  \begin{displaymath}
    A = P_r \cdot B \cdot P_c
  \end{displaymath}

\end{definition}

The lemma and corollary below provide us with basic results on
$\RELofTYPE{P}$ equivalence of \FKP\ products. A similar result can also
be found in \cite{VanLoan92}:
\bigskip

 
\begin{lemma}
\label{lem_four_split}

  Let $p$ and $q$ be relatively prime natural
  numbers. Then there exist permutation matrices $P_r,P_c$ such that
  \begin{equation}
    P_r \cdot  (F_p \otimes F_q) \cdot P_c = F_{pq}    \Longleftrightarrow
    (F_p \otimes F_q) = {P_r}^{-1} \cdot F_{pq} \cdot {P_c}^{-1} 
  \end{equation}

\end{lemma}



\PROOFstart
Let $\tilde{i}',\tilde{i}'';\tilde{j}',\tilde{j}''$ denote a {\bf shifted}
kronecker multiindex into a matrix with the kronecker product structure $p
\times q$, such as, for example, $F_p \otimes F_q$. The term {\bf
  shifted} here means that numbering starts from $0$. The relation
between shifted ordinary index $\tilde{i},\tilde{j}$
and shifted kronecker multiindex is, according to Definition \ref{def_kron_multiindex}:
\begin{equation}
\tilde{i}',\tilde{i}'';\tilde{j}',\tilde{j}''  \longleftrightarrow
\tilde{i}' \cdot q \; + \;  \tilde{i}'',\ \  \tilde{j}' \cdot q \; + \; \tilde{j}''  
\end{equation}
where
\begin{eqnarray}
  \tilde{i}' ,\tilde{j}'   \in  \{0 \ldots (p-1)\} \\
  \tilde{i}'',\tilde{j}''  \in  \{0 \ldots (q-1)\}
\end{eqnarray}


Now let the permutation matrices ${P_r}^{-1},{P_c}^{-1}$ move the
$\tilde{i}$-th row and $\tilde{j}$-th column, respectively, into the
specified below the $\tilde{i}',\tilde{i}''$-th row and
$\tilde{j}',\tilde{j}''$-th column of the result:
\begin{eqnarray}
\tilde{i}  \longrightarrow &  \tilde{i}',\tilde{i}''  =
\ \  (a\tilde{i}\mod p),\ \    (b \tilde{i} \mod  q)  \label{map_i_to_ii} \\
\tilde{j}  \longrightarrow &  \tilde{j}',\tilde{j}''  =  
\ \  (c\tilde{j}\mod p),\ \    (d \tilde{j} \mod  q)  \label{map_j_to_jj}
\end{eqnarray}
where $a,b,c,d$ are natural numbers satisfying, for some integers
$e,f$:
\begin{equation}
ep + (ac)q = 1 \ \ \ \ \ \  (bd)p + fq = 1   \label{eqty_rel_prime}
\end{equation}


Note that since $p,q$ are relatively prime, there exist integers
$e,x>0,y>0,f$ such that $ep + xq = 1$ and $yp + fq = 1$. To have $ac =
x$ and $bd = y$ we can take $a=1,\ c=x,\ b=1,\ d=y$.

Note also, that (\ref{eqty_rel_prime}) implies that $a,c$ are relatively
prime to $p$, and $b,d$ are relatively prime to $q$. Otherwise $1$
would have to have a divisor greater than $1$.\\


To show that maps (\ref{map_i_to_ii}),(\ref{map_j_to_jj}) properly define
permutations, that is that they are bijective, let us take
${\tilde{i}}_{1} \neq {\tilde{i}}_{2}$ belonging to $\{0, \ldots, (pq-1) \}$. If
they are mapped by (\ref{map_i_to_ii}) into equal pairs
$\tilde{i}',\tilde{i}''$, then:
\begin{eqnarray*}
  a{\tilde{i}}_{1} \mod p = a{\tilde{i}}_{2} \mod p  & \Longleftrightarrow
  &a({\tilde{i}}_{1}-{\tilde{i}}_{2}) \ \  \mbox{is divided by p} \\
  b{\tilde{i}}_{1} \mod q = b{\tilde{i}}_{2} \mod q  & \Longleftrightarrow
  &b({\tilde{i}}_{1}-{\tilde{i}}_{2}) \ \  \mbox{is divided by q}
\end{eqnarray*}
so ${\tilde{i}}_{1}-{\tilde{i}}_{2}$ is divided by $pq$ because pairs
$a,p$ and $b,q$ and $p,q$ are relatively prime. But this can't be so for
${\tilde{i}}_{1},{\tilde{i}}_{2} \in \{0, \ldots, (pq-1)\}$. Thus the map
(\ref{map_i_to_ii}) must be injective and similarily we show this for
the second map (\ref{map_j_to_jj}). Thus both maps map the sets $\{0,
  \ldots, (pq-1)\}$ and $\{0, \ldots, (p-1)\} \times \{0, \ldots,
  (q-1)\}$ one to one, that is they are bijective.\\


To show that $(F_p \otimes F_q) = {P_r}^{-1} \cdot F_{pq} \cdot
{P_c}^{-1}$ we will show that the $\tilde{i},\tilde{j}$-th element of
${F}_{pq}$ is mapped (or equal to) the
$\tilde{i}',\tilde{i}'';\tilde{j}',\tilde{j}''$-th, as defined in (\ref{map_i_to_ii})
and (\ref{map_j_to_jj}), element of $F_p
\otimes F_q$, i.e.:
\begin{equation}
\frac{1}{\sqrt{pq}} \exp{\left({\bf i} \frac{2\pi}{pq} \tilde{i} \tilde{j}\right)}
=
\frac{1}{\sqrt{p}} \exp{\left({\bf i} \frac{2\pi}{p} \tilde{i}' \tilde{j}'\right)}
\frac{1}{\sqrt{q}} \exp{\left({\bf i} \frac{2\pi}{q} \tilde{i}'' \tilde{j}''\right)}
\end{equation}
which is equivalent to the equality of phases:
\begin{equation}
  \frac{2\pi}{pq} \tilde{i} \tilde{j}   
     \stackrel{\mod 2\pi}{=}  
  \frac{2\pi}{pq} ( q\, \tilde{i}' \tilde{j}' \  + \  p\, \tilde{i}'' \tilde{j}'' ) 
\end{equation}
or to the equality for indices:  
\begin{equation} 
  \tilde{i} \tilde{j}
     \stackrel{\mod pq}{=}
  q\, \tilde{i}' \tilde{j}' \  + \  p\, \tilde{i}'' \tilde{j}''  \label{eqty_ind_mod_pq}
\end{equation}

{
\newcommand{\Ai}{a\,\tilde{i}}
\newcommand{\Bi}{b\,\tilde{i}}
\newcommand{\Cj}{c\,\tilde{j}}
\newcommand{\Dj}{d\,\tilde{j}}

\newcommand{\AiMp}{\Ai \mod p}
\newcommand{\BiMq}{\Bi \mod q}
\newcommand{\CjMp}{\Cj \mod p}
\newcommand{\DjMq}{\Dj \mod q}

\newcommand{\pijbis}{p\,\tilde{i}''\tilde{j}''}
\newcommand{\qijprim}{q\,\tilde{i}'\tilde{j}'}

Consider the difference:
\begin{eqnarray*}
\lefteqn{\tilde{i}\tilde{j} - \qijprim - \pijbis =}  \\
 & = & \frac{1}{ac} (\Ai)(\Cj)  - q(\AiMp)(\CjMp) - \pijbis   \\
 & = & \frac{1}{ac} (p\alpha + (\AiMp))(p\beta + (\CjMp)) - q(\AiMp)(\CjMp) - \pijbis \\  
 & = & \frac{1}{ac} (pA + (\AiMp)(\CjMp)) - q(\AiMp)(\CjMp) - \pijbis  \\
 & = & \frac{1}{ac} (pA + (\AiMp)(\CjMp)(1-qac)) - \pijbis  
\end{eqnarray*}
which is divided by $p$, for $1-qac = ep$ and $a,c$ are relatively
prime to $p$.

In the same way we show that it is divided by $q$:
\begin{eqnarray*}
\lefteqn{\tilde{i}\tilde{j} - \qijprim - \pijbis =}  \\
 & = & \frac{1}{bd} (\Bi)(\Dj) - \qijprim - p(\BiMq)(\DjMq) \\
 & = & \frac{1}{bd} (q\gamma + (\BiMq))(q\delta + (\DjMq)) - \qijprim - p(\BiMq)(\DjMq) \\  
 & = & \frac{1}{bd} (qB + (\BiMq)(\DjMq)) - \qijprim - p(\BiMq)(\DjMq)  \\
 & = & \frac{1}{bd} (qB + (\BiMq)(\DjMq)(1-pbd)) - \qijprim
\end{eqnarray*}
where $1-pbd = fq$ and $b,d$ are relatively prime to $q$.

Thus the considered difference is divided by relatively prime numbers
$p$ and $q$, so it is divided by $pq$, which is equivalent to (\ref{eqty_ind_mod_pq}).\PROOFend
}
\bigskip


\begin{corollary}
\label{col_four_split}

  Let $F$ be an \FKP\ product:
  \begin{equation}
    F\ \   =\ \   F_{n_{1}} \otimes F_{n_{2}} \otimes \ldots \otimes F_{n_{r}}
  \end{equation}
  Then there exist permutation matrices $P_r,P_c$ such that we have
  $\RELofTYPE{P}$ equivalence of matrices:
  \begin{equation}
    P_r \cdot F \cdot P_c  \ \  = \ \   F_{m_{1}} \otimes F_{m_{2}} \otimes \ldots \otimes F_{m_{s}}
  \end{equation}
  where the sequence $m_1, \ldots, m_s$ is obtained from the sequence
  $n_1, \ldots, n_r$ using a series of operations from the list below:
  \begin{enumerate}
   \item permuting a sequence, for example $n_1,n_2,n_3 \rightarrow
      n_1,n_3,n_2$ \label{permuting}
   \item merging a sequence: a subsequence $n_a,n_b$ can be replaced by
      $n_c = n_a n_b$ if $n_a,n_b$ are relatively prime \label{merging}
   \item division in a sequence: a sequence element $n_c$ can be replaced
      by a subsequence $n_a,n_b$ if $n_c = n_a n_b$, and $n_a,n_b$ are
      relatively prime \label{division}
  \end{enumerate}

\end{corollary}


\PROOFstart
Permuting the factors of a kronecker product, corresponding to operation
\ref{permuting} from the list, is equivalent to left and
right permuting this product (see Remark \ref{rem_reorder_kronp}).

Operations \ref{merging}, \ref{division} from the list correspond
to left and right permuting subproducts $F_{n_a} \otimes F_{n_b}$,
$F_{n_c}$ respectively, by Lemma \ref{lem_four_split}. Left and right permuting a
subproduct means left and right permuting the whole product, for
example
\begin{eqnarray}
A \otimes (P_1 B P_2) \otimes C \ & = &\  A \otimes (P_1 \otimes I_C)(B\otimes C)(P_2 \otimes I_C) \nonumber \\
 & = & (I_A \otimes P_1 \otimes I_C)(A \otimes B \otimes C)(I_A \otimes P_2\otimes I_C) \nonumber
\end{eqnarray}
where $A,B,C$ square matrices, $I_A,I_B,I_C$ the identity matrices of
the same size, respectively.

Combining all the permutation matrices corresponding to the operations
performed on the size sequence (factor sequence) leads to $P_r,P_c$.\PROOFend
\bigskip


\begin{example}

  Some $\RELofTYPE{P}$-equivalence sets of \FKP\ products.

\end{example}


For an example we will split the sets of all \FKP\ products of size 16
and 72 into $\RELofTYPE{P}$ equivalence subsets, with their all
elements featured up to the order of factors, using the allowed
operations of Corollary \ref{col_four_split}.

However, we do not show at this moment that matrices from different sets below
are not $\RELofTYPE{P}$ equivalent,  i.e. that the sets below are
$\RELofTYPE{P}$ equivalence classes.
\medskip


The case of $n_1 \cdot \ldots \cdot n_r \ \ \   = \ \ \  16\ (=2 \cdot 2 \cdot
2 \cdot 2)$:
\begin{enumerate}
\item $F_2 \otimes F_2 \otimes F_2 \otimes F_2$

\item $F_4 \otimes F_2 \otimes F_2$

\item $F_4 \otimes F_4$

\item $F_8 \otimes F_2$

\item $F_{16}$
\end{enumerate}


The case of $n_1 \cdot \ldots \cdot n_r \ \ \  = \ \ \  72\ (=2 \cdot 2 \cdot
2 \cdot 3 \cdot 3)$:
\begin{enumerate}
\item $F_2 \otimes F_2 \otimes F_2 \otimes F_3 \otimes F_3$, \ \  $F_2
      \otimes F_2 \otimes F_6 \otimes F_3$, \ \  $F_2 \otimes F_6
      \otimes F_6$

\item $F_2 \otimes F_2 \otimes F_2 \otimes F_9$, \ \  $F_2 \otimes F_2
      \otimes F_{18}$

\item $F_2 \otimes F_4 \otimes F_3 \otimes F_3$, \ \  $F_2 \otimes
      F_{12} \otimes F_3$, \ \  $F_6 \otimes F_4 \otimes F_3$, \ \
      $F_6 \otimes F_{12}$

\item $F_2 \otimes F_4 \otimes F_9$, \ \  $F_2 \otimes F_{36}$, \ \
      $F_4 \otimes F_{18}$

\item $F_8 \otimes F_3 \otimes F_3$, \ \  $F_{24} \otimes F_3$

\item $F_8 \otimes F_9$, \ \  $F_{72}$
\end{enumerate}
\bigskip

We continue with pure ordered \FFKP's, see Lemma
\ref{lem_poFFKP_row_types}.
\medskip


\begin{lemma}
\label{lem_type_row_number_change}

  Let there be an $m$ factor pure ordered \FFKP\ given:
  \begin{equation}
    F = F_{a^{k_m}} \otimes \ldots \otimes F_{a^{k_1}}
  \end{equation}
  with $\NofROWS{p}{(m)}{a}{0}$ rows of type $\ROWtype{a}{0}$, $\NofROWS{p}{(m)}{a}{1}$ rows
  of type $\ROWtype{a}{1}$, ..., $\NofROWS{p}{(m)}{a}{k_m}$ rows of type
  $\ROWtype{a}{k_m}$.

  Then a new pure ordered \FFKP\ product:
  \begin{equation}
    F' = F_{a^{k_{m+1}}} \otimes F,\ \ \ k_{m+1} \geq k_m
  \end{equation}
  contains rows of types $\ROWtype{a}{0}, \ROWtype{a}{1}, \ldots,
  \ROWtype{a}{k_{m+1}}$ in quantities given by the formulas:
  \begin{eqnarray}
    \NofROWS{p}{(m+1)}{a}{0} & = & 1 \cdot \NofROWS{p}{(m)}{a}{0} \\
    \NofROWS{p}{(m+1)}{a}{k}  & = & a^{k-1} (a-1) \left(
      \NofROWS{p}{(m)}{a}{0} + \NofROWS{p}{(m)}{a}{1} + \ldots +
      \NofROWS{p}{(m)}{a}{k-1} 
    \right) \  +\   a^k \cdot \NofROWS{p}{(m)}{a}{k}
    \label{form_already_present_count}\\
    & & \mbox{for $1 < k \leq k_m$}   \nonumber \\
    \NofROWS{p}{(m+1)}{a}{l}  & = & a^{l-1} (a-1) \left(
      \NofROWS{p}{(m)}{a}{0} + \NofROWS{p}{(m)}{a}{1} + \ldots +
      \NofROWS{p}{(m)}{a}{k_m} 
    \right) \label{form_newly_introduced_count} \\
    & = &  a^{l-1}(a-1)a^{\sum_{i=1}^{m} k_i} \nonumber \\
    & & \mbox{for $k_m < l \leq k_{m+1}$}  \nonumber
  \end{eqnarray}

\end{lemma}


\PROOFstart
From Lemma \ref{lem_poFFKP_row_types} $F$ contains rows of types
$\ROWtype{a}{0}, \ROWtype{a}{1}, \ldots, \ROWtype{a}{k_m}$ and $F'$ has
rows of types $\ROWtype{a}{0}, \ROWtype{a}{1}, \ldots,
\ROWtype{a}{k_{m+1}}$. 

We calculate the number of rows of type $\ROWtype{a}{k}$,
where $1< k \leq k_m$, in $F'$ using Lemma
\ref{lem_row_kronp_type}. Let $\NofROWS{p}{}{a}{n}$ denote the number of
type $\ROWtype{a}{n}$ rows in $F_{a^{k_{m+1}}}$, which is given by
Lemma \ref{lem_four_rows_count}, and let $\mathcal{P}$ denote the
cartesian product $\{0,1,\ldots,k_{m+1}\} \times
\{0,1,\ldots,k_{m}\}$. Then:
\begin{eqnarray*}
  \lefteqn{\NofROWS{p}{(m+1)}{a}{k} =} \\ 
  & = & \sum_{\{(g,h) \in \mathcal{P}:\max{(g,h)}=k\}} \NofROWS{p}{}{a}{g}
         \cdot  \NofROWS{p}{(m)}{a}{h} \\
  & = & \sum_{h=0}^{k-1} \NofROWS{p}{}{a}{k} \cdot \NofROWS{p}{(m)}{a}{h} \ 
      + \ \sum_{g=0}^{k} \NofROWS{p}{}{a}{g} \cdot \NofROWS{p}{(m)}{a}{k} \\
  & = & a^{k-1} (a-1) \left( \NofROWS{p}{(m)}{a}{0} + \NofROWS{p}{(m)}{a}{1} + \ldots +
        \NofROWS{p}{(m)}{a}{k-1} \right) + \\
  &   & + (1 + (a-1) + \ldots + (a-1)a^{k-1}) \cdot \NofROWS{p}{(m)}{a}{k} \\
  & = & a^{k-1} (a-1) \left( \NofROWS{p}{(m)}{a}{0} + \NofROWS{p}{(m)}{a}{1} + \ldots +
        \NofROWS{p}{(m)}{a}{k-1} \right) + a^{k} \cdot \NofROWS{p}{(m)}{a}{k}
\end{eqnarray*}

Rows of type $\ROWtype{a}{0}$ are obtained in $F'$ only as kronecker
products of rows of type $\ROWtype{a}{0}$, so from Lemma
\ref{lem_four_rows_count}:
\begin{displaymath}
  \NofROWS{p}{(m+1)}{a}{0} = \NofROWS{p}{}{a}{0} \cdot
  \NofROWS{p}{(m)}{a}{0} = 1 \cdot \NofROWS{p}{(m)}{a}{0}
\end{displaymath}

The number of rows of type $\ROWtype{a}{l}$, where $k_m < l \leq
k_{m+1}$, if there are any in $F'$ (if $k_{m+1} > k_m$), is
calculated in the following way, where $\mathcal{P}$ defined above:
\begin{eqnarray*}
  \lefteqn{\NofROWS{p}{(m+1)}{a}{l} =} \\ 
  & = & \sum_{\{(g,h) \in \mathcal{P}:\max{(g,h)}=l\}} \NofROWS{p}{}{a}{g}
       \cdot  \NofROWS{p}{(m)}{a}{h} \\
  & = & \sum_{h=0}^{k_m} \NofROWS{p}{}{a}{l} \cdot \NofROWS{p}{(m)}{a}{h} \\
  & = & a^{l-1} (a-1) \left( \NofROWS{p}{(m)}{a}{0} + \NofROWS{p}{(m)}{a}{1} + \ldots +
      \NofROWS{p}{(m)}{a}{k_m} \right)  \\
  & = &  a^{l-1}(a-1)a^{(k_1+\ldots+k_m)}
\end{eqnarray*}

This completes the proof.\PROOFend
\bigskip

We will also need the definition of an \DEFINED{introduction index}:
\bigskip


\begin{definition}
\label{def_intro_index}

  An \DEFINED{introduction index} for type $\ROWtype{a}{k}$ rows in $F$ being an ordered
  pure \FFKP\ product of $F_{a^n}$ matrices, denoted $\INTROix{a}{k}^F$, is the
  position (number) of the factor of $F$, counting from left to right
  starting from $1$ for the first factor, which first introduces rows of type
  $\ROWtype{a}{k}$ in the process of constructing $F$ by left
  kronecker multiplying by consecutive factors of nondecreasing size.

\end{definition}
\bigskip


\begin{example}

  Introduction indices for $F_{16} \otimes F_{8} \otimes F_{8} \otimes
  F_{2}$

\end{example}


Take the ordered pure \FFKP\ : $F = F_{16} \otimes F_{8} \otimes F_{8} \otimes
  F_{2}$. Then its introduction indices are:
\begin{equation}
\INTROix{2}{0}^F = 1,\ \ 
\INTROix{2}{1}^F = 1,\ \ 
\INTROix{2}{2}^F = 2,\ \ 
\INTROix{2}{3}^F = 2,\ \ 
\INTROix{2}{4}^F = 4
\end{equation}
\bigskip


\begin{remark}
  Let $F$ and $F'=F_{a^{k_{m+1}}} \otimes F$ be ordered pure \FFKP\
  products of $F_{a^n}$ matrices. Then introduction indices for $F'$ are
  inherited from $F$. 
\end{remark}

Because of the remark above, in the proofs below we will usually ommit
the upper index $F$ when writing introduction indices for a p.o. \FFKP
\ $F$, as we
will often have to consider p.o. \FFKP\ subproducts of $F$.
\bigskip


\begin{theorem}
\label{theor_poFFKP_type_row_count}

  Let $F$ be a pure ordered \FFKP, consisting of $s$ factors:
  \begin{equation}
    F = F_{a^{k_s}} \otimes \ldots \otimes F_{a^{k_1}}
  \end{equation}

  Then $F$ has $\NofROWS{p}{(s)}{a}{0}$ rows of type $\ROWtype{a}{0}$,
  $\NofROWS{p}{(s)}{a}{1}$ rows of type $\ROWtype{a}{1}$, ...,
  $\NofROWS{p}{(s)}{a}{m}$ rows of type $\ROWtype{a}{m}$, ..., where
  $0 \leq m \leq k_s$, and the numbers $\NofROWS{p}{(s)}{a}{0},
  \NofROWS{p}{(s)}{a}{1},\ldots,\NofROWS{p}{(s)}{a}{m},\ldots$ are
  given by:
  \begin{eqnarray*}
    \NofROWS{p}{(s)}{a}{0} & = & 1          \\
    \NofROWS{p}{(s)}{a}{1} & = & a^s -1      \\
    \NofROWS{p}{(s)}{a}{2} & = & a^{2s - (\Nt{2}-1)} - a^s  \\
    \NofROWS{p}{(s)}{a}{3} & = & a^{3s - (\Nt{3}-1) -
      (\Nt{2}-1)}  -   a^{2s - (\Nt{2}-1)}   \\
    \NofROWS{p}{(s)}{a}{4} & = & a^{4s - (\Nt{4}-1) -
      (\Nt{3}-1) - (\Nt{2}-1)}  -   a^{3s - (\Nt{3}-1) -
      (\Nt{2}-1)}   \\
    & \ldots &  \\
    \NofROWS{p}{(s)}{a}{m} & = & a^{ms - \sum_{l=0}^{m} (\Nt{l}-1)}  -  
                             a^{(m-1)s - \sum_{l=0}^{m-1} (\Nt{l}-1)} \\
    & \ldots &
  \end{eqnarray*}
  where $\Nt{0}=1$, $\Nt{1}=1$, $1 \leq \Nt{2} \leq \ldots \leq \Nt{m}
  \leq \ldots \leq \Nt{k_s}$ are the introduction indices of rows of
  the type $\ROWtype{a}{0}, \ROWtype{a}{1}, \ROWtype{a}{2}, \ldots,
  \ROWtype{a}{m}, \ldots, \ROWtype{a}{k_s}$, respectively.

  Each $m$-th formula giving
  $\NofROWS{p}{(s)}{a}{m}$, where $1 \leq m \leq k_s$, is valid for $s
  \geq (\Nt{m}-1)$, the threshhold argument value used when applying
  the formula to the right $(\Nt{m}-1)$ factor subproduct of $F$.

\end{theorem}




\PROOFstart
The type $\ROWtype{a}{0}$ and $\ROWtype{a}{1}$ rows are present in any
$F_{a^n}$, that is why we put \mbox{$\Nt{0}=1,\ \Nt{1}=1$,} i.e. these row
types are introduced by the first right factor of a pure ordered \FFKP\
of $F_{a^n}$ matrices.


Lemma \ref{lem_type_row_number_change} gives us the rules for changing
$\NofROWS{p}{(s)}{a}{0}$ when adding left new factors to a
p.o. \FFKP\:
\begin{equation}
  \NofROWS{p}{(s+1)}{a}{0} = 1 \cdot \NofROWS{p}{(s)}{a}{0}
\end{equation}
 Because in the first, Fourier matrix, factor we have $1$
row of type $\ROWtype{a}{0}$ (see Lemma \ref{lem_four_rows_count}), so we
have $\NofROWS{p}{(s)}{a}{0}=1$ for any $s$ and any $s$-factor
p.o.\FFKP\ of $F_{a^n}$ matrices, for any $a$ prime.


The number of type $\ROWtype{a}{1}$ rows in the first factor is equal to
$\NofROWS{p}{(1)}{a}{1}=(a-1)$ (again Lemma
\ref{lem_four_rows_count}), and using Lemma
\ref{lem_type_row_number_change} we get the rule:
\begin{equation}
  \NofROWS{p}{(s+1)}{a}{1} = (a-1)\NofROWS{p}{(s)}{a}{0} + a\NofROWS{p}{(s)}{a}{1}
\end{equation}
which yields $\NofROWS{p}{(s)}{a}{1} = a^s - 1$ for any $s$-factor
p.o.\FFKP\ of $F_{a^n}$ matrices. Since $\NofROWS{p}{(0)}{a}{1} = 0$
for the $0$-factor subproduct, we regard $\NofROWS{p}{(s)}{a}{1}$
as valid for $s \geq 0 = \Nt{1}-1$.


Now assume that the formulas for $\NofROWS{p}{(s)}{a}{k}$ are correct
for $k=0,1,2,\ldots,m$ and that they are valid for $s \geq (\Nt{k}-1)$,
respectively. We will show that for $s \geq (\Nt{m+1}-1)$ and $m \geq 1$:
\begin{equation}
\label{form_next_type_rows_count}
  \NofROWS{p}{(s)}{a}{m+1} = a^{(m+1)s - \sum_{l=0}^{m+1} (\Nt{l}-1)}  -  
                             a^{ms - \sum_{l=0}^{m} (\Nt{l}-1)}
\end{equation}
Note that even $\NofROWS{p}{(s)}{a}{1}$ can be written using the
general scheme:
\begin{equation}
\label{form_type_a_count_general}
  \NofROWS{p}{(s)}{a}{1} = a^s - 1 = a^{1s - (\Nt{1}-1) - (\Nt{0}-1)}
  - a^{0s - (\Nt{0}-1)}
\end{equation}



{\bf Case} $s = \Nt{m+1} - 1$

We start from calculating $\NofROWS{p}{(s)}{a}{m+1}$ for $s = \Nt{m+1}
- 1$ substituted into \ref{form_next_type_rows_count}. The result is
$0$ as it should for there are no rows of type $\ROWtype{a}{m+1}$
before they are introduced by the $\Nt{m+1}$-th factor of a p.o.\FFKP.
It is $0$ also if $\Nt{m+1}=1$ (we can interprete
$\NofROWS{p}{(0)}{a}{m+1}$ as the number of type $\ROWtype{a}{m+1}$
rows in no product at all), so $\NofROWS{p}{(0)}{a}{m+1}$ if used in a
sum (fortunately it won't be) will not cause any damage because of
abuse of index $(s)$.
\bigskip


{\bf Case} $s = \Nt{m+1}$

To check whether the formula (\ref{form_next_type_rows_count}) gives
the correct value for $s=\Nt{m+1}$ we have to consider two subcases.


Firstly, let $s = \Nt{m+1} = \Nt{m} = \Nt{m-1} = \ldots = \Nt{0} = 1$,
which means that rows of type $\ROWtype{a}{m+1}$ are introduced by the
first factor of a p.o.\FFKP\ and from Lemma \ref{lem_four_rows_count}
their number is $a^m (a-1)$. Let us substitute all those $1$'s into
(\ref{form_next_type_rows_count}). We will get the required value:
\begin{equation}
  a^{(m+1) \cdot 1} - a^{m \cdot 1} = a^m (a-1)
\end{equation}


Secondly, let $s = \Nt{m+1} > 1$. Recall that we are talking about the
situation in which the so far constructed ($s$ factors) p.o. \FFKP\ has
just had type $\ROWtype{a}{m+1}$ rows introduced by the $\Nt{m+1}$-th
($s$-th) left factor. That is why we have to use the formula
(\ref{form_newly_introduced_count}) from
Lemma \ref{lem_type_row_number_change} to calculate the resulting
number of type $\ROWtype{a}{m+1}$ rows. Note that we do not know 'the
highest' present row type, so we write (which will be shown to be
correct in a moment):
\begin{equation}
\label{form_next_introduced_count}
  \NofROWS{p}{(\Nt{m+1})}{a}{m+1} = a^m (a-1) \left(
    \NofROWS{p}{(\Nt{m+1}-1)}{a}{0} + \NofROWS{p}{(\Nt{m+1}-1)}{a}{1}
    + \ldots + \NofROWS{p}{(\Nt{m+1}-1)}{a}{m} \right)
\end{equation}

Note that even if for some right $\NofROWS{p}{(\Nt{m+1}-1)}{a}{k}$
(except for the first two left) we have that $\Nt{m} = \ldots =
\Nt{m-d}$ are all equal to $\Nt{m+1}$, that is the types
$\ROWtype{a}{m+1},\ROWtype{a}{m},\ldots,\ROWtype{a}{m-d}$ are
introduced together, the corresponding \\
$\NofROWS{p}{(\Nt{m+1}-1)}{a}{m},\ \ldots,\ 
\NofROWS{p}{(\Nt{m+1}-1)}{a}{m-d}$ values are in fact the values of
$\NofROWS{p}{(\Nt{m}-1)}{a}{m},\ \ldots,\ 
\NofROWS{p}{(\Nt{m-d}-1)}{a}{m-d}$. 
From the induction assumption on
correct formulas for $\NofROWS{p}{(s)}{a}{k}$,\ $s \geq
\Nt{k}-1$, $k=0,\ldots,m$ we can apply safely these formulas to calculate
$\NofROWS{p}{(\Nt{m+1})}{a}{m+1}$ from
(\ref{form_next_introduced_count}). The values
$\NofROWS{p}{(\Nt{m}-1)}{a}{m},\ \ldots,\ 
\NofROWS{p}{(\Nt{m-d}-1)}{a}{m-d}$ will in this case be zeros as they
should.

The long sum in (\ref{form_next_introduced_count}), from the induction
assumption formulas, is equal to
\begin{equation}
  a^{m(\Nt{m+1}-1) - \sum_{l=0}^{m} (\Nt{l}-1)}
\end{equation}
Note that it is also true if $m=1$ at the begining of the induction
process, as (\ref{form_type_a_count_general}) is satisfied.

So, what on one hand we get from (\ref{form_next_introduced_count}) is
\begin{equation}
  \NofROWS{p}{(\Nt{m+1})}{a}{m+1} = a^m (a-1) a^{m(\Nt{m+1}-1) - \sum_{l=0}^{m} (\Nt{l}-1)}
\end{equation}
and on the other, substituting $\Nt{m+1}$ into
(\ref{form_next_type_rows_count}), we obtain the same:
\begin{eqnarray*}
  \lefteqn{\NofROWS{p}{(\Nt{m+1})}{a}{m+1} =}\\ 
   & = & a^{(m+1)\Nt{m+1} - \sum_{l=0}^{m+1} (\Nt{l}-1)}  -  
         a^{m\Nt{m+1} - \sum_{l=0}^{m} (\Nt{l}-1)} \\
   & = & a^{m(\Nt{m+1}-1) - \sum_{l=0}^{m} (\Nt{l}-1)} 
          \left( a^{\Nt{m+1} + m - (\Nt{m+1}-1)} - a^m \right) \\
   & = & a^{m(\Nt{m+1}-1) - \sum_{l=0}^{m} (\Nt{l}-1)} (a^m (a-1))
\end{eqnarray*}


{\bf Case} $s > \Nt{m+1}$

Consider the formula (\ref{form_next_type_rows_count}) for $s >
\Nt{m+1}$, that is type $\ROWtype{a}{m+1}$ is introduced by a factor
preceding the $s$-th one. Since formula
(\ref{form_next_type_rows_count}) is correct, as we have shown above,
for $s = \Nt{m+1}$, we can use formula
(\ref{form_already_present_count}) from Lemma
\ref{lem_type_row_number_change} to obtain next values of
$\NofROWS{p}{(s)}{a}{m+1}$ for $s = \Nt{m+1}+1,\ \ \Nt{m+1}+2,\ \ 
\ldots$. We will show that in each step of this process we get
$\NofROWS{p}{(s)}{a}{m+1}$ of the form
(\ref{form_next_type_rows_count}). We start from
$\NofROWS{p}{(s)}{a}{m+1}$ given by (\ref{form_next_type_rows_count}),
then from (\ref{form_already_present_count}) we get:
\begin{equation}
  \NofROWS{p}{(s+1)}{a}{m+1} = 
  a^m (a-1) \left( \NofROWS{p}{(s)}{a}{0} + \NofROWS{p}{(s)}{a}{1} +
    \ldots + \NofROWS{p}{(s)}{a}{m} \right)
  + a^{m+1} \NofROWS{p}{(s)}{a}{m+1}
\end{equation}
where the long sum in brackets can be safely replaced, even for $m=1$,
by $a^{ms - \sum_{l=0}^{m} (\Nt{l}-1)}$ from the main induction assumption,
then $\NofROWS{p}{(s+1)}{a}{m+1}$ is given by:
\begin{eqnarray*}
  \lefteqn{\NofROWS{p}{(s+1)}{a}{m+1} =}  \\
  & = & a^m (a-1) a^{ms - \sum_{l=0}^{m} (\Nt{l}-1)} +
        a^{m+1} \left( a^{(m+1)s - \sum_{l=0}^{m+1} (\Nt{l}-1)} -
        a^{ms - \sum_{l=0}^{m} (\Nt{l}-1)}         \right) \\
  & = & a^{(m+1)(s+1) - \sum_{l=0}^{m+1} (\Nt{l}-1)} -
        a^{m(s+1) - \sum_{l=0}^{m} (\Nt{l}-1)}
\end{eqnarray*}
which is analogous to (\ref{form_next_type_rows_count}), with $s$
replaced by $s+1$. This, combined with
(\ref{form_next_type_rows_count}) working for $s = \Nt{m+1}$, means
that (\ref{form_next_type_rows_count})
works also for $s > \Nt{m+1}$.
\medskip


We have thus shown that (\ref{form_next_type_rows_count}) gives
correct values of the number of rows of type $\ROWtype{a}{m+1}$ for $s
\geq \Nt{m+1} - 1$, which completes the proof by induction.\PROOFend
\bigskip


\begin{theorem}
\label{theor_poFFKP_identity_criteria}
  Let there be two pure ordered \FFKP\ products of $F_{a^n}$ matrices
  given ($a$ prime):
  \begin{eqnarray*}
    F' = & F_{a^{k_{s_1}}} \otimes F_{a^{k_{s_1 - 1}}} \otimes \ldots \otimes
    F_{a^{k_{1}}} &
    \ \ \ k_{s_1} \geq k_{s_1 - 1} \geq \ldots \geq k_{1} \\
    F'' = & F_{a^{l_{s_2}}} \otimes F_{a^{l_{s_2 - 1}}} \otimes \ldots \otimes
    F_{a^{l_{1}}} &
    \ \ \ l_{s_2} \geq l_{s_2 - 1} \geq \ldots \geq l_{1}
  \end{eqnarray*}

  Then the following statements \ITEMa, \ITEMb,
  \ITEMc\ are equivalent:

  \begin{description}

    \item[\ITEMa]
          Both \FFKP\ products, $F'$ and $F''$, contain:
          \begin{itemize}
            \item equal numbers of type $\ROWtype{a}{0}$ rows:
                  $\NofROWS{p'}{(s_1)}{a}{0} =
                  \NofROWS{p''}{(s_2)}{a}{0}$
            \item equal numbers of type $\ROWtype{a}{1}$ rows:
                  $\NofROWS{p'}{(s_1)}{a}{1} =
                  \NofROWS{p''}{(s_2)}{a}{1}$
            \item $\ldots$
            \item equal numbers of type $\ROWtype{a}{k_{s_1}} =
                  \ROWtype{a}{l_{s_2}}$ rows:
                  $\NofROWS{p'}{(s_1)}{a}{k_{s_1}} =
                  \NofROWS{p''}{(s_2)}{a}{l_{s_2}}$
           \end{itemize}

    \item[\ITEMb]
          Introduction indices $\Nt{0}, \Nt{1}, \ldots, \Nt{k_{s_1} =
            l_{s_2}}$ are common to both $F'$ and $F''$, and $s_1 =   
          s_2$.

    \item[\ITEMc]
          The p.o. \FFKP\ products $F'$ and $F''$ are identical, that
          is $s_1 = s_2$ and $k_1 = l_1$, $k_2 = l_2$, ..., $k_{s_1} =
          l_{s_2}$.
  \end{description}

\end{theorem}


\PROOFstart
We consider the cases of implication:
\smallskip

  
{\bf Case } \ITEMa\ $\Longrightarrow$ \ITEMb.
 
Introduction indices $\Nt{0}=1$ and $\Nt{1}=1$ are common to both $F'$
and $F''$, for rows of type $\ROWtype{a}{0}$ and $\ROWtype{a}{1}$ are
introduced by any first right factor $F_{a^n},\  n \geq 1$.

Assume that common to $F'$ and $F''$ are the introduction indices
$\Nt{0}, \Nt{1}, \ldots, \Nt{k} \in \{1,2,\ldots,\min{(s_1,s_2)}\}$, where $k < k_{s_1} =
l_{s_2}$. Note that the last equality results from the statement
\ITEMa\ and the last inequality must be in accordance with Lemma
\ref{lem_poFFKP_row_types}.  We will show now
that $\Nt{k+1}$ is also common to $F'$ and $F''$, that is $\Nt{k+1}' =
\Nt{k+1}''$.
\smallskip

Rows of type $\ROWtype{a}{k+1}$ are introduced into $F'$ by the $\Nt{k+1}'$-th
factor of $F'$ and into $F''$ by $\Nt{k+1}''$-th factor of $F''$,
where $\Nt{k+1}' \leq s_1$, $\Nt{k+1}'' \leq s_2$ and $k+1 \leq
k_{s_1} = l_{s_2}$. Theorem \ref{theor_poFFKP_type_row_count} provides
us with the written below formulas for the number of type
$\ROWtype{a}{k+1}$ 
rows in $F'$ and $F''$. From the statement \ITEMa\ these
numbers are equal, which is expressed by:
\begin{eqnarray*}
  \lefteqn{a^{(k+1){s_1} - (\Nt{k+1}'-1) - \sum_{l=0}^{k} (\Nt{l}-1)}  - 
  a^{ k{s_1} - \sum_{l=0}^{k} (\Nt{l}-1)} =}  \\
  & & a^{(k+1){s_2} - (\Nt{k+1}''-1) - \sum_{l=0}^{k} (\Nt{l}-1)}  - 
  a^{ k{s_2} - \sum_{l=0}^{k} (\Nt{l}-1)}
\end{eqnarray*}
which is equivalent to
\begin{equation}
\label{eq_equal_row_count_products}
  a^{k{s_1} - \sum_{l=0}^{k} (\Nt{l}-1)}
  \left( a^{{s_1} - (\Nt{k+1}'-1)} - 1 \right)
  \ \ =\ \ 
  a^{k{s_2} - \sum_{l=0}^{k} (\Nt{l}-1)}
  \left( a^{{s_2} - (\Nt{k+1}''-1)} - 1 \right)
\end{equation}

From the uniqueness of prime number factorization of a natural number
and the fact that $a$ is prime, powers of $a$ on the left and right
side of (\ref{eq_equal_row_count_products}) are equal, so $s_1 = s_2$,
which we get in each step of the induction process. As a result the
factors not divided by $a$ in (\ref{eq_equal_row_count_products}) are
also equal, what is more $s_1 = s_2$, so $\Nt{k+1}' = \Nt{k+1}'' =
\Nt{k+1}$.

Thus all $\Nt{k},\  k \in \{0,1,2,\ldots,k_{s_1}=l_{s_2}\}$ are common
to $F'$ and $F''$, 
and the implication \ITEMa\ $\Longrightarrow$ \ITEMb\
is proved.
\medskip

{\bf Case } \ITEMb\ $\Longrightarrow$ \ITEMc.

From the statement \ITEMb\ we have that $s_1 = s_2$.

Now, assume that $F'$ and $F''$ are not identical. That is, for some
$s \leq s_1 = s_2$ the \mbox{$1$-st}, \mbox{$2$-nd}, ..., \mbox{$(s-1)$-th} right factors of
$F'$ and $F''$ are equal and the \mbox{$s$-th} factors are not. Say, they are
$F_{a^k}$ for $F'$ and $F_{a^l}$ for $F''$ such that $k > l$. Then the
introduction indices $\Nt{k}',\Nt{k}''$ for $F',F''$ must satisfy
$\Nt{k}' \leq s < \Nt{k}''$. The reason is that in $F''$ rows of type
$\ROWtype{a}{k}$ have not yet been introduced by the \mbox{$s$-th} factor,
whereas in $F'$ they have already been introduced by some \mbox{$(s-t)$-th}
factor, $t \geq 0$. So $F'$ and $F''$ cannot have all their
introduction indices common.
\medskip


{\bf Case } \ITEMc\ $\Longrightarrow$ \ITEMb. 

Obvious.
\medskip


{\bf Case } \ITEMb\ $\Longrightarrow$ \ITEMa.

The implication is given by\ \  
\mbox{\ITEMb\ $\Longrightarrow$ \ITEMc\
  $\Longrightarrow$ \ITEMa},\ \  where\ \  
\mbox{\ITEMc\ $\Longrightarrow$ \ITEMa}\ \  is obvious.
\medskip


{\bf Cases } \mbox{\ITEMa\ $\Longrightarrow$ \ITEMc}\ \ 
and\ \  \mbox{\ITEMc\ $\Longrightarrow$ \ITEMa}

It is obvious that\ \  \mbox{\ITEMc\ $\Longrightarrow$
  \ITEMa}.\ \  Further,\ \  \mbox{\ITEMa\ $\Longrightarrow$
  \ITEMb}\ \  and\ \  \mbox{\ITEMb\ $\Longrightarrow$
  \ITEMc}\ \  implies\ \  \mbox{\ITEMa\ $\Longrightarrow$
  \ITEMc}.

\PROOFend
\medskip


\begin{corollary}
  Pure ordered \FFKP\ products $F'$ and $F''$ are $\RELofTYPE{P}$
  equivalent if and only if they are identical ($\Longleftrightarrow$
  equal).
\end{corollary}


\PROOFstart
If $F'$ and $F''$ are $\RELofTYPE{P}$ equivalent, they satisfy the
statement \ITEMa\ of Theorem
\ref{theor_poFFKP_identity_criteria}, which, by this theorem, is
equivalent to $F'$ anf $F''$ being identical ($\Leftrightarrow$ equal
by Theorem \ref{theor_poFFKP_identity_criteria}).

The opposite implication is obvious.\PROOFend
\bigskip


\begin{theorem}
\label{theor_FFKP_identity_criteria}
  Let there be two \FFKP\ products $F'$ and $F''$ given (see
  Definition \ref{def_FFKP}) of size $N = a_1^{b_1} \cdot a_2^{b_2} \cdot
  \ldots \cdot a_r^{b_r}$, where $a_1 > a_2 > \ldots > a_r$ are prime
  factors of $N$. 

  Then $F'$ and $F''$ are $\RELofTYPE{P}$ equivalent
  if and only if $F'$ and $F''$ are identical up to the order of their
  factors ($\Longleftrightarrow$ their corresponding maximal pure
  ordered \FFKP\ subproducts are identical in the sense of statement
  \ITEMc\ of Theorem \ref{theor_poFFKP_identity_criteria}).
\end{theorem}


{
\newcommand{\Sone}[1]{s_1^{(#1)}}
\newcommand{\Stwo}[1]{s_2^{(#1)}}
\newcommand{\gexp}[2]{g_{#1}^{(#2)}}
\newcommand{\hexp}[2]{h_{#1}^{(#2)}}

\PROOFstart
It is enough to prove Theorem \ref{theor_FFKP_identity_criteria} in
the case of $F'$ and $F''$ ordered (which by Remark
\ref{rem_reorder_kronp} corresponds to left and
right permuting) in such a way that:
\begin{eqnarray*}
  F' & = & F^{\prime(1)} \otimes F^{\prime(2)} \otimes \ldots \otimes F^{\prime(r)} \\
  F''& = & F^{\bis(1)} \otimes F^{\bis(2)} \otimes \ldots \otimes F^{\bis(r)} \\
\end{eqnarray*}
where $F^{\prime(k)}$ and $F^{\bis(k)}$ are pure ordered \FFKP's such that
\begin{eqnarray*}
  F^{\prime(k)} & = & F_{a_k^{\gexp{\Sone{k}}{k}}}    \otimes
                      F_{a_k^{\gexp{\Sone{k}-1}{k}}}  \otimes
                                          \ldots  \otimes
                      F_{a_k^{\gexp{1}{k}}}  \\
  F^{\bis(k)} & = & F_{a_k^{\hexp{\Stwo{k}}{k}}}    \otimes
                    F_{a_k^{\hexp{\Stwo{k}-1}{k}}}  \otimes
                                        \ldots  \otimes
                    F_{a_k^{\hexp{1}{k}}}  \\
\end{eqnarray*}
where
\begin{equation}
  \sum_{s=1}^{\Sone{k}} \gexp{s}{k} \  =\ \sum_{s=1}^{\Stwo{k}}
  \hexp{s}{k} \ =\  b_k 
\end{equation}

Let us denote by
$\NofROWS{p'}{(\Sone{k})}{a_k}{m},\ \NofROWS{p''}{(\Stwo{k})}{a_k}{m}$
the numbers of type 
$\ROWtype{a_k}{m}$ rows in 
$\Sone{k}$ factor p.o. \FFKP\ subproduct $F^{\prime(k)}$ and 
$\Stwo{k}$ factor p.o. \FFKP\ subproduct $F^{\bis(k)}$, respectively.

From Lemma \ref{lem_rel_prime_type_row_kronp_type} 
each row of type $\ROWtype{a_k}{m}$ of $F'$ 
is obtained as a kronecker product of 
type $\ROWtype{a_l}{0}$ rows 
from $F^{\prime(l)}$ where $l \in \{1,2,\ldots,r\}\setminus\{k\}$ 
and a row of type $\ROWtype{a_k}{m}$
from $F^{\prime(k)}$.
Similarly for $F''$. So, the numbers of type $\ROWtype{a_k}{m}$ rows
in $F'$ and $F''$ will be, by Theorem
\ref{theor_poFFKP_type_row_count}:
\begin{eqnarray*}
  \left( \prod_{l \in \{1,2,\ldots,r\} \setminus \{k\}}
              \NofROWS{p'}{(\Sone{l})}{a_l}{0} \right)
  \cdot \NofROWS{p'}{(\Sone{k})}{a_k}{m} 
                     & = &
                       \NofROWS{p'}{(\Sone{k})}{a_k}{m}   \\
  \left( \prod_{l \in \{1,2,\ldots,r\} \setminus \{k\}}
              \NofROWS{p''}{(\Stwo{l})}{a_l}{0} \right)
  \cdot \NofROWS{p''}{(\Stwo{k})}{a_k}{m}
                     & = &
                       \NofROWS{p''}{(\Stwo{k})}{a_k}{m}  \\     
\end{eqnarray*}

If $F'$ and $F''$ are $\RELofTYPE{P}$ equivalent, then for each $k \in
\{1,2,\ldots,r\}$
the above numbers of type $\ROWtype{a_k}{m}$ rows in $F'$ and $F''$
are equal, where $m \in
\{0,1,\ldots,\gexp{\Sone{k}}{k}=\hexp{\Stwo{k}}{k}\}$ (the last equality in
the set description follows from $\RELofTYPE{P}$ equivalence of
$F',F''$ and Lemmas \ref{lem_rel_prime_type_row_kronp_type},
\ref{lem_poFFKP_row_types}). That is the pairs
$\NofROWS{p'}{(\Sone{k})}{a_k}{m}, \NofROWS{p''}{(\Stwo{k})}{a_k}{m}$ are
pairs of equal numbers, $k$ and $m$ as above. Then, from Theorem
\ref{theor_poFFKP_identity_criteria} $F^{\prime(k)}$ and $F^{\bis(k)}$ are
identical (equal) for each $k \in
\{1,2,\ldots,r\}$. This means that $F'$ and $F''$, assumed ordered,
are identical.
\smallskip

If $F'$ and $F''$ are identical, they are $\RELofTYPE{P}$ equivalent
of course.
\smallskip

Now, taking Remark \ref{rem_reorder_kronp} into consideration, we can
generalize the above implications to unordered \FFKP's, to obtain Theorem
\ref{theor_FFKP_identity_criteria}.\PROOFend
}
\bigskip

%
%

\section{Permutation-phasing equivalence classes of kronecker products of Fourier matrices}
\label{sec_fkp_PD_equivalence}

Our next step will be generalizing Theorem
\ref{theor_FFKP_identity_criteria} to the case of permutation-phasing
equivalence, denoted by $\RELofEQUI$ and defined below. This is an
equivalence relation.


\begin{definition}
\label{def_permutation_phasing_equivalence}
  Two square matrices $A$ and $B$ of the same size are \DEFINED{permutation-phasing
    equivalent}, i.e. $A \RELofEQUI B$ if and only if there exist
  permutation matrices $P_r,\ P_c$ and unitary diagonal matrices
  $D_r,\ D_c$ such that
  \begin{displaymath}
    A = P_r \cdot D_r \cdot B \cdot D_c \cdot P_c
  \end{displaymath}

  Note that if $A \RELofEQUI B$ then
  \begin{displaymath}
    A = 
    (P_r D_r P_r^T) \cdot P_r \cdot B \cdot P_c \cdot (P_c^T D_c P_c)
    = D_r' \cdot P_r \cdot B \cdot P_c \cdot D_c' 
  \end{displaymath}
  where $D_r',\ D_c'$ are unitary diagonal matrices.
\end{definition}

We will prove
that left and right multiplying an \FKP\ (\FFKP\ as a special case) of
size $N$ by such unitary diagonal matrices that the result still has a
row and column filled with $1/\sqrt{N}$ is equivalent to permuting this
\FKP. We start with a single Fourier matrix, using the lemma below.
\medskip


\begin{lemma}
\label{lem_DrDc_uniqueness}
  
  Let $U$ be a unitary matrix of size $n$ and let $D_r,\ D_c$ be unitary diagonal
  matrices such that $[D_c]_{1,1} = 1$ and the elements of the $i$-th
  row and $j$-th column of the matrix:
  \begin{equation}
    U' = D_r U D_c
  \end{equation}
  have prescribed phases.

  Then $D_r,\ D_c$ are uniquely determined.
\end{lemma}


\PROOFstart 
Let $\alpha_1,\ldots,\alpha_n$ denote the prescribed phases in the
$i$-th row of $U'$ and $\beta_1,\ldots,\beta_n$ the prescribed phases
in the $j$-th column of $U'$, $\alpha_j = \beta_i$. Let also
$\gamma_{k,l}$ denote the phase of $U_{k,l}$, $\phi_{k}$ the phase of
$[D_r]_{k,k}$ and $\psi_{l}$ the phase of $[D_c]_{l,l}$.

Since $[D_c]_{1,1} = 1 \Rightarrow \psi_{1} = 0$, there must be
$\gamma_{i,1} + \phi_i \stackrel{\mod 2\pi}{=} \alpha_1$, which
determines $\phi_i$ up to multiplicity of $2\pi$.
This, in turn, determines the phases $\psi_2,\ldots,\psi_n$ of $D_c$
as they must satisfy $\gamma_{i,l} + \phi_i + \psi_l \stackrel{\mod
  2\pi}{=} \alpha_l$. Among others, $\psi_j$ is thus determined. If
$j=1$, $\psi_j = 0$ from the assumption.

Finally, the phases $\phi_1,\ldots,\phi_{i-1},\phi_{i+1},\ldots,\phi_n$
are forced to satisfy $\gamma_{k,j} + \phi_k + \psi_j \stackrel{\mod
  2\pi}{=} \beta_k$, where $k \in \{1,\ldots,n\} \setminus \{ i \}$,
while $\gamma_{i,j} + \phi_i + \psi_j \stackrel{\mod 2\pi}{=} \beta_i
= \alpha_j$
is satisfied as a result of the proper choice of $\psi_j$ in the
previous step.\PROOFend
\bigskip


\begin{lemma}
\label{lem_Fourier_phasing_is_shifting}
  
  Let $F_n$ be a Fourier matrix of size $n$ and let $D_r,\ D_c$ be
  unitary diagonal matrices such that the matrix $F_n'$:
  \begin{equation}
    F_n' = D_r F_n D_c
  \end{equation}
  has its $i$-th row and $j$-th column filled with
  $\frac{1}{\sqrt{n}}$ (phases equal to $0$).

  Then
  \begin{equation}
    F_n' = P_r F_n P_c
  \end{equation}
  where $P_r,\ P_c$ are permutation matrices such that:
  \begin{description}
    \item[] $P_r$ moves rows indexed by $1,2,\ldots,n$ into positions\ \
          \mbox{$i,i+1,\ldots,n,1,2,\ldots,i-1$}
    \item[] $P_r$ moves columns indexed by $1,2,\ldots,n$ into positions\ \
          \mbox{$j,j+1,\ldots,n,1,2,\ldots,j-1$}
  \end{description}

\end{lemma}


\PROOFstart
We can assume that $[D_c]_{1,1} = 1$. Otherwise, for $D_r,\ D_c$ there
exist $D'_r,\ D'_c$ such that $[D'_c]_{1,1}=1$ and $D'_r A D'_c = D_r A
D_c$ for any matrix $A$, and we consider $D'_r,\ D'_c$ instead of
$D_r,\ D_c$.

We will use shifted indices (from $0$ to $n-1$) from now on, written with
tildes.
\smallskip

Let us take a ``shifted'' Fourier matrix, denoted by $F_n''$:
\begin{equation}
\label{def_Fbis_F_shifted}
  [F_n'']_{\Si{k},\Si{l}} 
          = [F_n]_{\Si{k} - \Si{i} \mod n,\Si{l} - \Si{j} \mod n}
          = \frac{1}{\sqrt{n}} e^{\IfAf{n} (\Si{k} - \Si{i}
            \mod n)(\Si{l} - \Si{j} \mod n)}
\end{equation}
where $\Si{i} = i-1$, $\Si{j} = j-1$.

We will show that it can be obtained by left and right multiplying
$F_n$ by unitary diagonal matrices $D''_r,\ D''_c$ such that
$[D''_c]_{1,1}=1$. Namely, we will prove that:
\begin{equation}
\label{eq_Fbis_F_phased}
  [F_n'']_{\Si{k},\Si{l}} = [F_n]_{\Si{k},\Si{l}} \cdot
                            e^{-\IfAf{n} (\Si{k}-\Si{i})\Si{j}} \cdot
                            e^{-\IfAf{n} \Si{l} \Si{i}}
                          =
                            \frac{1}{\sqrt{n}} 
                            e^{\IfAf{n} (\Si{k}\Si{l} -
                              (\Si{k}-\Si{i})\Si{j} - \Si{l}\Si{i})}
\end{equation} 

To do it we need to prove the equality modulo $n$ of expressions from
the phases of (\ref{def_Fbis_F_shifted}) and (\ref{eq_Fbis_F_phased}):
\begin{equation}
  (\Si{k} - \Si{i}\mod n)(\Si{l} - \Si{j} \mod n)
  \stackrel{\mod n}{=}
  \Si{k}\Si{l} - (\Si{k}-\Si{i})\Si{j} - \Si{l}\Si{i}  
\end{equation}
The above holds since the left side is equal modulo $n$ to $(\Si{k} -
\Si{i})(\Si{l} - \Si{j})$, which in turn is equal to the right
side.

Thus $F_n''$ is given by:
\begin{equation}
  D''_r \cdot F_n \cdot D''_c \ \ \ \mbox{where}\ 
  [D''_r]_{\Si{k},\Si{k}} = e^{-\IfAf{n} (\Si{k}-\Si{i})\Si{j}},\ 
  [D''_c]_{\Si{l},\Si{l}} = e^{-\IfAf{n} \Si{l} \Si{i}},\ 
  [D''_c]_{\Si{0},\Si{0}} = 1
\end{equation}
and $F_n''$ has its $i$-th rows and $j$-th column filled with
$1/\sqrt{n}$, where $i,j$ are normal (not shifted) indices.

$F_n''$ satisfies the assumption imposed on $F_n'$ in Lemma
\ref{lem_Fourier_phasing_is_shifting}, then by \mbox{Lemma
\ref{lem_DrDc_uniqueness}}\ \ $D_r = D''_r$, $D_c = D''_c$ are determined
uniquely, so $F_n' = F_n''$. The last equality is equivalent to the
statement of Lemma \ref{lem_Fourier_phasing_is_shifting}.\PROOFend
\bigskip

We generalize the above result to the set of all \FKP\ products:


\begin{lemma}
\label{lem_FKP_phasing_is_permuting}

  Let $F$ be an \FKP\ product of the form:
  \begin{equation}
    F = F_{N_1} \otimes \ldots \otimes F_{N_p}
  \end{equation}
  and let $D_r,\ D_c$ be such unitary diagonal matrices that the
  matrix:
  \begin{equation}
    F' = D_r \cdot F \cdot D_c
  \end{equation}
  has its $i$-th row and $j$-th column filled with $1/{\sqrt{N}}$,
  where $N = N_1 \cdot \ldots \cdot N_p$ is the size of $F$.

  Then there exist permutation matrices $P_r,\ P_c$ such that
  \begin{equation}
    F' = P_r \cdot F \cdot P_c
  \end{equation}

\end{lemma}


\PROOFstart
We assume that $[D_c]_{1,1} = 1$. Otherwise we take $D'_r,\ D'_c$ such
that $[D'_c]_{1,1} = 1$ and for any matrix $A$ there holds $D'_r A
D'_c = D_r A D_c$, instead of $D_r,\ D_c$.

From all the assumptions on $D_r,\ D_c$, by Lemma
\ref{lem_DrDc_uniqueness} the matrices $D_r,\ D_c$ are unique.

We will construct the matrices $D_r,\ D_c$ satisfying the above
assumptions.
\smallskip

Let the kronecker multiindex $i_1,\ldots,i_p;j_1,\ldots,j_p$ 
correspond to the index pair $i,\ j$ of 
the selected row and column of $F$, see Definition \ref{def_kron_multiindex}. 
Now, let unitary diagonal matrix pairs $D^{(k)}_r,\ D^{(k)}_c$,
$k=1,\ldots,p$ be such that $[D^{(k)}_c]_{1,1} = 1$ and the matrix
\begin{equation}
  F'_{N_k} = D^{(k)}_r \cdot F_{N_k} \cdot D^{(k)}_c
\end{equation}
has its $i_k$-th row and $j_k$-th column filled with $1/\sqrt{N_k}$,
for $k = 1,\ldots,p$. By Lemma \ref{lem_DrDc_uniqueness} the matrices
$D^{(k)}_r,\ D^{(k)}_c$ are unique, and they are given by Lemma 
\ref{lem_Fourier_phasing_is_shifting}. 

As the matrix
\begin{eqnarray*}
 \lefteqn{F'_{N_1} \otimes \ldots \otimes F'_{N_p} =} \\
 & = &     (D^{(1)}_r F_{N_1} D^{(1)}_c) \otimes \ldots \otimes
            (D^{(p)}_r F_{N_p} D^{(p)}_c)  \\
 & = &  (D^{(1)}_r \otimes \ldots \otimes D^{(p)}_r)
         (F_{N_1}   \otimes \ldots \otimes F_{N_p}  )
         (D^{(1)}_c \otimes \ldots \otimes D^{(p)}_c)
\end{eqnarray*}
has its $i$-th row and $j$-th column ($i,j$ corresponding bijectively to the
kronecker multiindex $i_1,\ldots,i_p;j_1,\ldots,j_p$) filled with
$1/\sqrt{N_1 \cdot \ldots \cdot N_p} = 1/\sqrt{N}$, then by Lemma
\ref{lem_DrDc_uniqueness} we have that
\begin{eqnarray*}
  D_r & = & (D^{(1)}_r \otimes \ldots \otimes D^{(p)}_r) \\
  D_c & = & (D^{(1)}_c \otimes \ldots \otimes D^{(p)}_c)
\end{eqnarray*}
and
\begin{equation}
\label{eq_Fprim_kronprod_of_phased_factors}
  F'= (F'_{N_1}   \otimes \ldots \otimes F'_{N_p})
\end{equation}

On the other hand, by Lemma \ref{lem_Fourier_phasing_is_shifting}, for
$k = 1,\ldots,p$ there exist permutation matrices $P^{(k)}_r,\
P^{(k)}_c$ such that
\begin{equation}
  F'_{N_k} = P^{(k)}_r \cdot F_{N_k} \cdot P^{(k)}_c
\end{equation}

So, from (\ref{eq_Fprim_kronprod_of_phased_factors}):
\begin{eqnarray*}
 \lefteqn{F' =} \\
 & = & (P^{(1)}_r F_{N_1} P^{(1)}_c) \otimes \ldots \otimes
       (P^{(p)}_r F_{N_p} P^{(p)}_c) \\
 & = & (P^{(1)}_r \otimes \ldots \otimes P^{(p)}_r)
       (F_{N_1}   \otimes \ldots \otimes F_{N_p}  )
       (P^{(1)}_c \otimes \ldots \otimes P^{(p)}_c) \\
 & = & P_r \cdot F \cdot P_c
\end{eqnarray*}
which completes the proof.\PROOFend
\bigskip

Using the above two Lemmas \ref{lem_Fourier_phasing_is_shifting} and
\ref{lem_FKP_phasing_is_permuting} we can prove that:


\begin{theorem}
\label{theor_FKPs_PD_equiv_is_P_equiv}

  Two \FKP\ products $F'$ and $F''$ are $\RELofTYPE{P}$ equivalent if
  and only if they are $\RELofEQUI$ equivalent.

\end{theorem}



\PROOFstart
If $F'$ and $F''$ are $\RELofTYPE{P}$ equivalent, then obviously they
are $\RELofEQUI$ equivalent. We take identity matrices for $D_r,\ D_c$
in Definition \ref{def_permutation_phasing_equivalence}. 
\smallskip

If $F'$ and $F''$ are $\RELofEQUI$ equivalent, then for
some permutation matrices $P_r,\ P_c$ and unitary diagonal matrices
$D_r,\ D_c$ :
\begin{equation}
\label{eq_Fprim_Fbis_PD_related}
  F'' = P_r \cdot D_r F' D_c \cdot P_c
\end{equation}

$F''$, as an \FKP\ product, has the first row and column filled with
$1/\sqrt{N}$, where $N$ is the size of $F'$ and $F''$, so the matrix
$D_r F' D_c$ in (\ref{eq_Fprim_Fbis_PD_related}) must have the $i$-th
row and $j$-th column filled with $1/\sqrt{N}$ as well, for some $i,j
\in \{1,\ldots,N\}$. Then by Lemma \ref{lem_FKP_phasing_is_permuting}
there exist permutation matrices $P'_r,\ P'_c$ such that:
\begin{equation}
  D_r F' D_c = P'_r F' P'_c
\end{equation}.

From the above and (\ref{eq_Fprim_Fbis_PD_related}) $F''$ is equal to
\begin{equation}
  (P_r P'_r) \cdot F' \cdot (P'_c P_c)
\end{equation}
which means that $F'$ and $F''$ are $\RELofTYPE{P}$
equivalent.\PROOFend
\bigskip
\bigskip

Let $\FKPs$ denote the set of all \FKP\ products and let
$\FKPs_N$ denote the set of all \FKP\ products of size $N$.

By Theorem \ref{theor_FKPs_PD_equiv_is_P_equiv} the $\RELofTYPE{P}$
equivalence classes within $\FKPs$ are the same as $\RELofEQUI$
equivalence classes in $\FKPs$. Further, by Corollary \ref{col_four_split} (on operations
on \FKP 's preserving $\RELofTYPE{P}$ equivalence) and by Theorem
\ref{theor_FFKP_identity_criteria}  (on the 
criteria of two \FFKP 's being $\RELofTYPE{P}$ equivalent),
these classes are represented by ordered \FFKP\ products with ordered
pure \FFKP\ subproducts, that is by \FFKP's of the form:
\begin{equation}
  \left( F_{a_1^{k_{m^{(1)}}^{(1)}}} \otimes 
                              \ldots \otimes
         F_{a_1^{k_{1}^{(1)}}}               \right)  \otimes 
                                              \ldots  \otimes
  \left( F_{a_r^{k_{m^{(r)}}^{(r)}}} \otimes 
                              \ldots \otimes
         F_{a_r^{k_{1}^{(r)}}}               \right)
\end{equation}
where $a_1,a_2,\ldots,a_r$ are prime numbers satisfying
\begin{equation}
  a_1 > a_2 > \ldots > a_r
\end{equation}
and $k_{m}^{(l)}$,\ \ $m = 1,\ldots,m^{(l)}$,\ \ $l = 1,\ldots,r$\ \ are
respective exponents such that
\begin{equation}
  k_{m^{(l)}}^{(l)} \geq k_{m^{(l)}-1}^{(l)} \geq \ldots \geq k_{1}^{(l)}\
  \ \ \mbox{for  $l = 1,\ldots,r$}
\end{equation}

To calculate the number of $\RELofEQUI$ equivalence classes of \FKP 's of a given
size $N$, that is in $\FKPs_N$, we need the definition of a partition:


\begin{definition}
\label{def_partition}

  A \DEFINED{partition} of a natural number $N$ is any ordered
  sequence of positive ($>0$) integers $(n_1,n_2,\ldots,n_p)$
  such that\ \  \mbox{$n_1 \geq n_2 \geq \ldots \geq n_p$}\ \  and\ \  \mbox{$n_1
    + n_2 + \ldots + n_p = N$}.

\end{definition}

The number of partitions of $N$ is usually denoted by $\NofPART{N}$
and we will use this notation. A good reference on partitions is
\cite{Partitions}.

From what was said above about the representatives of $\RELofEQUI$
equivalence classes in $\FKPs$ we immediately have:


\begin{theorem}

  The number of $\RELofEQUI$ equivalence classes within $\FKPs_N$
  (\FKP's of size $N$), where $N$ has the prime number factorization:
  \begin{equation}
    N = a_1^{b_1} \cdot a_2^{b_2} \cdot \ldots \cdot a_r^{b_r},
    \ \ \ \ \ \
    a_1 > a_2 > \ldots > a_r
  \end{equation}
  is equal to the product of the numbers of partitions:
  \begin{equation}
    \NofPART{b_1} \cdot \NofPART{b_2} \cdot \ldots \cdot \NofPART{b_r}
  \end{equation}

\end{theorem}
\smallskip

We will denote the number of $\RELofEQUI$ equivalence classes in
$\FKPs_N$ by $\NofEQcl{N}$. Below we present a few examples of sets of
all $\RELofEQUI$ equivalence classes for the values of $N$: $30, 48,
36$. $\EQclass{}{F}$ denotes the class represented by $F$. For each
class, all (up to the order of factors) \FKP's contained in it are listed.


\begin{example}
  Sets of $\RELofEQUI,\ (\RELofTYPE{P})$ equivalence classes within
  $\FKPs_N$,\ \ \  $N = 30,\ 48,\ 36$.
\end{example}


\begin{description}

  \item[$\FKPs_{5 \cdot 3 \cdot 2}\ :$]         
      \begin{description}
        \item[]  $\EQclass{}{F_5 \otimes F_3 \otimes F_2}$ contains:
                     \begin{displaymath}
                       F_5 \otimes F_3 \otimes F_2,\ \ \ \ 
                       F_6 \otimes F_5,\ \ \ \ 
                       F_{10} \otimes F_3,\ \ \ \ 
                       F_{15} \otimes F_2,\ \ \ \ 
                       F_{30}
                     \end{displaymath}
      \end{description}

  \item[$\FKPs_{3 \cdot 2^4}\ :$]
      \begin{description}   
        \item[]  $\EQclass{}{F_3 \otimes F_{2^4}}$ contains:
                     \begin{displaymath}
                        F_{16} \otimes F_3,\ \ \ \ 
                        F_{48}
                     \end{displaymath}

        \item[]  $\EQclass{}{F_3 \otimes F_{2^3} \otimes F_2}$ contains:
                     \begin{displaymath}
                       F_8 \otimes F_2 \otimes F_3,\ \ \ \ 
                       F_{24} \otimes F_2,\ \ \ \ 
                       F_8 \otimes F_6
                     \end{displaymath}

        \item[]  $\EQclass{}{F_3 \otimes F_{2^2} \otimes F_{2^2}}$ contains:
                     \begin{displaymath}
                       F_4 \otimes F_4 \otimes F_3,\ \ \ \ 
                       F_{12} \otimes F_4,\ \ \ \ 
                     \end{displaymath}

        \item[]  $\EQclass{}{F_3 \otimes F_{2^2} \otimes F_2 \otimes F_2}$ contains:
                     \begin{displaymath}
                       F_4 \otimes F_2 \otimes F_2 \otimes F_3,\ \ \ \ 
                       F_{12} \otimes F_2 \otimes F_2,\ \ \ \ 
                       F_4 \otimes F_2 \otimes F_6
                     \end{displaymath}

        \item[]  $\EQclass{}{F_3 \otimes F_2 \otimes F_2 \otimes F_2 \otimes F_2}$ contains:
                     \begin{displaymath}
                       F_3 \otimes F_2 \otimes F_2 \otimes F_2 \otimes F_2,\ \ \ \ 
                       F_6 \otimes F_2 \otimes F_2 \otimes F_2
                     \end{displaymath}
      \end{description}

  \item[$\FKPs_{3^2 \cdot 2^2}\ :$]
      \begin{description} 
        \item[]  $\EQclass{}{F_{3^2} \otimes F_{2^2}}$ contains:
                     \begin{displaymath}
                       F_9 \otimes F_4,\ \ \ \ 
                       F_{36}
                     \end{displaymath}

        \item[]  $\EQclass{}{F_{3^2} \otimes F_2 \otimes F_2}$ contains:
                     \begin{displaymath}
                       F_9 \otimes F_2 \otimes F_2,\ \ \ \ 
                       F_{18} \otimes F_2
                     \end{displaymath}

        \item[]  $\EQclass{}{F_3 \otimes F_3 \otimes F_{2^2}}$ contains:
                     \begin{displaymath}
                       F_4 \otimes F_3 \otimes F_3,\ \ \ \ 
                       F_{12} \otimes F_3
                     \end{displaymath}

        \item[]  $\EQclass{}{F_3 \otimes F_3 \otimes F_2 \otimes F_2}$ contains:
                     \begin{displaymath}
                       F_3 \otimes F_3 \otimes F_2 \otimes F_2,\ \ \ \ 
                       F_6 \otimes F_3 \otimes F_2,\ \ \ \ 
                       F_6 \otimes F_6
                     \end{displaymath}
      \end{description}

\end{description}

%
%

\section*{Acknoledgements}
\label{sec_acknoledgements}

This work was supported by the Polish Ministry of Scientific Research
and Information Technology under the (solicited) grant No PBZ-Min-008/P03/2003.

%
%

\end{document}